# Explicit and Exact Traveling Wave Solutions of Cahn Allen equation using MSE Method

[1]Harun-Or-Roshid , [2]M. Zulfikar Ali and [1]Md. Rafiqul Islam
[1]Department of Mathematics, Pabna University of Science and Technology, Bangladesh
[2]Department of Mathematics, Rajshahi University, Rajshahi, Bangladesh
[*]Email: harunorroshidmd@gmail.com or harun_math@pust.ac.bd

**Abstract.** By using Modified simple equation method, we study the Cahn Allen equation which arises in many scientific applications such as mathematical biology, quantum mechanics and plasma physics. As a result, the existence of solitary wave solutions of the Cahn Allen equation is obtained. Exact explicit solutions interms of hyperbolic solutions of the associated Cahn Allen equation are characterized with some free parameters. Finally, the variety of structure and graphical representation make the dynamics of the equations visible and provides the mathematical foundation in mathematical biology, quantum mechanics and plasma physics.

## 1. Introduction

The mathematical modeling of events in nature can be explained by differential equations. It is well familiar that various types of the physical phenomena in the field fluid mechanics, quantum mechanics, electricity, plasma physics, chemical kinematics, propagation of shallow water waves and optical fibers are modeled by nonlinear evolution equation and the appearance of solitary wave solutions in nature is somewhat frequent. But, the nonlinear processes are one of the major challenges and not easy to control because the nonlinear characteristic of the system abruptly changes due to some small changes in parameters including time. Thus, the issue becomes more difficult and hence crucial solution is needed. The solutions of these equations have crucial impact in mathematical physics and engineering. The variety of solutions of NLEEs, that are mutual operating different mathematical techniques, is very important in many fields of science such as fluid mechanics, optical fibers, technology of space, control engineering problems, hydrodynamics, meteorology, plasma physics, applied mathematics. Advance nonlinear techniques are major to solve inherent nonlinear problems, particularly those are involving dynamical systems and allied areas. In recent years, there become large improvements in finding the exact solutions of NLEEs. Many powerful methods have been established and enhanced, such as, the modified extended Fan sub-equation method [1], the homogeneous balance method [2-3], the Jacobi elliptic function expansion [4], the Backlund transformation method [5, 6], the Darboux transformation method [7],the Adomian decomposition method [8-9], the auxiliary equation method[10, 11], the $(G'/G)$-expansion method [12-18], the Exp-$(-\varphi(\xi))$-expansion method [19], the sine-cosine method [20-22], the tanh method [23], the F-expansion method [24, 25], the exp-function method [26, 30], the modified simple equation method [27–29], first integral method [31], Simple equation method [32], Bilinear method [33], transformed rational function method [34] and so on. Most of the above methods are depend on computational software except MSE method.

The objective of this paper is to look for new exact traveling wave solutions including topological soliton, single soliton solutions of the well-recognized Cahn-Allen equation [30, 31] via MSE method.

## 2. Description of the MSE Method

Consider a general form of a nonlinear evolution equation,

$$H(u, u_t, u_x, u_{xt}, u_{xx}, \cdots \cdots) = 0. \tag{1}$$

where $u(\xi) = u(x,t)$ is an unknown function, $H$ is a polynomial of $u(x,t)$ and its partial derivatives in which the highest order derivatives and nonlinear terms are involved. In the following, we present the main steps of the method:

**Step 1:** Combine the real variables $x$ and $t$ by a compound variable $\xi$

$$u(\bar{r},t) = u(\xi), \quad \xi = \bar{P}.\bar{r} \pm wt, \tag{2}$$

Here $\bar{P} = l\hat{i} + m\hat{j} + n\hat{k}$ and $\bar{P} = x\hat{i} + y\hat{j} + z\hat{k}$ where $l, m, n$ are the constant magnitudes along the axes of $x, y, z$ respectively, $\mu$ is wave number and $w$ is the speed of the traveling wave.

This travelling wave transformation permits us to reduce Eq. (1) to the following ordinary differential equation (ODE):

$$G(u, u', u'' \cdots \cdots \cdots) = 0. \tag{3}$$

where $G$ is a polynomial in $u(\xi)$ and its derivatives, wherein $u'(\xi) = \dfrac{du}{d\xi}, u''(\xi) = \dfrac{d^2u}{d\xi^2}$ and so on.

**Step 2:** We also consider that the Eq.(3) has the formal solution

$$u(\xi) = \sum_{i=0}^{n} A_i \left( \frac{S'(\xi)}{S(\xi)} \right)^i. \tag{4}$$

where $A_i (0 \le i \le n)$ are constants to be determined, and $S(\xi)$ is also unknown function to be evaluated.

**Step 3:** The value of positive integer $n$ in Eq. (4) can be determined by taking into account the homogeneous balance between the highest order nonlinear terms and the derivatives of highest order occurring in Eq. (3). If the degree of $u(\xi)$ is, $D(u(\xi)) = n$ then the degree of the other expression will be $u(\xi)$ as follows:

$$D\left( \frac{d^p u(\xi)}{d\xi^p} \right) = n + p, \quad D\left( u^p \left( \frac{d^q u(\xi)}{d\xi^q} \right)^s \right) = np + s(n+q).$$

**Step 4**: Inserting Eq. (4) into Eq. (3), we get a polynomial of $(S'(\xi)/S(\xi))$ and its derivatives and $(S(\xi))^{-i}, (i = 0,1,2,\cdots,n)$. In the resultant polynomial, we equate all the coefficients of $(S(\xi))^{-i}, (i = 0,1,2,\cdots,n)$ to zero. This technique produces a system of algebraic and differential equations which can be solved receiving $A_i (i = 0,1,2,\cdots,n)$, $S(\xi)$ and the value of the other needful parameters. This completes the determination of the solution to the Eq. (1).

**Remark:** In comparison the modified simple equation method with the simple equation method [32], it is seen that Simple equation method take help of auxiliary equation (Riccati equation) but modified simple equation method can perform directly without help of auxiliary equation. On the other hand, Simple equation gives results which are special case of Modified equation method.

**3. Traveling wave solution of Cahn Allen Equation**

Let's consider nonlinear parabolic partial differential equation given by

$$u_t = u_{xx} - u^m + u. \tag{5}$$

for $m = 3$, Eq.(5) becomes Cahn Allen equation [28, 29]. This equation arises in many scientific applications such as mathematical biology, quantum mechanics and plasma physics. To solve this example, we can use transformation $\xi = kx + wt$ (where, $k$ and $w$ are the wave number and the wave speed, respectively) then Eq. (5) becomes to an ordinary differential equation

$$wu' - k^2 u'' + u^3 - u = 0. \tag{6}$$

Balancing $u^3$ with $u''$ then gives $n = 1$.

$$u(\xi) = A_0 + A_1 \frac{S'(\xi)}{S(\xi)}, \tag{7}$$

$$u'(\xi) = A_1 \frac{S''(\xi)}{S(\xi)} - A_1 \left(\frac{S'(\xi)}{S(\xi)}\right)^2, \tag{8}$$

$$u''(\xi) = A_1 \frac{S'''(\xi)}{S(\xi)} - 3A_1 \frac{S''(\xi)S'(\xi)}{S^2(\xi)} + 2A_1 \left(\frac{S'(\xi)}{S(\xi)}\right)^3. \tag{9}$$

Putting Eq.(6)-Eq.(9) in the Eq. (6) and equating coefficients of like powers of $\frac{S'(\xi)}{S(\xi)}$, we get

Coefficient of $(S(\xi))^0$: $\quad A_0^3 - A_0 = 0,$ \hfill (10)

Coefficient of $(S(\xi))^{-1}$: $-k^2 A_1 S'''(\xi) + 3A_0^2 A_1 S'(\xi) + wA_1 S''(\xi) - A_1 S'(\xi) = 0,$ \hfill (11)

Coefficient of $(S(\xi))^{-2}$: $-wA_1(S'(\xi))^2 + 3k^3 A_1 S'(\xi)S''(\xi) + 3A_0 A_1^2(S'(\xi))^2 = 0,$ \hfill (12)

Coefficient of $(S(\xi))^{-3}$: $A_1(A_1^2 - 2k^2)(S'(\xi))^3 = 0.$ \hfill (13)

From Eq.(10), we achieve $A_0 = 0, 1, -1$ and from Eq.(13), $A_1 \neq 0$ thus $A_1 = \pm\sqrt{2}k$,

$$\frac{S'''}{S''} = \frac{3k^2(3A_0^2 - 1) + w(w - 3A_0 A_1)}{k^2(w - 3A_0 A_1)}. \tag{14}$$

Integrating we have

$$S'' = c_1 \exp\left(\frac{3k^2(3A_0^2 - 1) + w(w - 3A_0 A_1)}{k^2(w - 3A_0 A_1)}\xi\right), \tag{15}$$

From 12 $S' = \frac{3c_1 k^2}{w - 3A_0 A_1} \exp\left(\frac{3k^2(3A_0^2 - 1) + w(w - 3A_0 A_1)}{k^2(w - 3A_0 A_1)}\xi\right),$ \hfill (16)

and $S = \frac{3c_1 k^4}{3k^2(3A_0^2 - 1) + w(w - 3A_0 A_1)} \exp\left(\frac{3k^2(3A_0^2 - 1) + w(w - 3A_0 A_1)}{k^2(w - 3A_0 A_1)}\xi\right) + c_2.$ \hfill (17)

Using (16) and (17), we attain

$$u = A_0 + \frac{3c_1 A_1 k^2}{w - 3A_0 A_1} \times \frac{\exp\left(\frac{3k^2(3A_0^2 - 1) + w(w - 3A_0 A_1)}{k^2(w - 3A_0 A_1)}\xi\right)}{\frac{3c_1 k^4}{3k^2(3A_0^2 - 1) + w(w - 3A_0 A_1)} \exp\left(\frac{3k^2(3A_0^2 - 1) + w(w - 3A_0 A_1)}{k^2(w - 3A_0 A_1)}\xi\right) + c_2} \tag{18}$$

where $\xi = k(x \pm \frac{3}{\sqrt{2}}t)$ with $w = \pm\frac{3}{\sqrt{2}}k$. Here $c_1$ and $c_2$ are arbitrary constants.

**Case-I:** For set $A_0 = 0$, $A_1 = \pm\sqrt{2}k$, we get

$$u = \pm \frac{3\sqrt{2}c_1 k^3}{w} \times \frac{\exp(\frac{w^2 - 3k^2}{k^2 w}\xi)}{\frac{3c_1 k^4}{w^2 - 3k^2}\exp(\frac{w^2 - 3k^2}{k^2 w}\xi) + c_2}, \tag{19}$$

where $\xi = k(x \pm \frac{3}{\sqrt{2}}t)$ with $w = \pm \frac{3}{\sqrt{2}}k$.

If $c_2 = \frac{3k^4 c_1}{w^2 - 3k^2}$ then $u = \pm \frac{w^2 - 3k^2}{\sqrt{2}wk}\left\{1 + \tanh\left(\frac{w^2 - 3k^2}{2wk^2}\xi\right)\right\}, \tag{20}$

where $\xi = k(x \pm \frac{3}{\sqrt{2}}t)$ with $w = \pm \frac{3}{\sqrt{2}}k$.

If $c_2 = -\frac{3k^4 c_1}{w^2 - 3k^2}$ then $u = \pm \frac{w^2 - 3k^2}{\sqrt{2}wk}\left\{1 + \coth\left(\frac{w^2 - 3k^2}{2wk^2}\xi\right)\right\}, \tag{21}$

where $\xi = k(x \pm \frac{3}{\sqrt{2}}t)$.

Since $c_1$ and $c_2$ are arbitrary constants, for other choices of $c_1$ and $c_2$ it might yield much new and more general exact solutions of the nonlinear Cahn Allen equation without any aid of symbolic computation software. The solutions u(x, t) obtained in Eqs. (20) and (21) are presented in the following figures: (see Figs. 1 and 2).

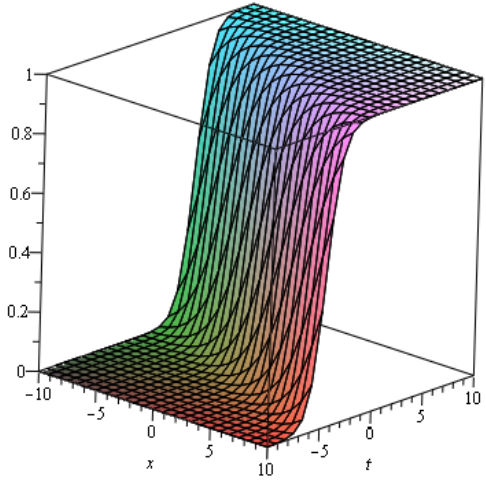 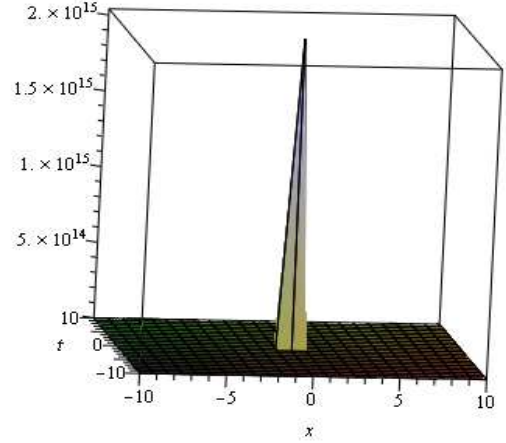

Fig-1: Kink wave of the solution Eq.(20) with $k=1$

Fig-2: Single soliton solution of the Eq.(21) with $k=1$

**Case-II:** For set $A_0 = \pm 1$, $A_1 = \pm\sqrt{2}k$, we get

$$u = \pm 1 \pm \frac{3\sqrt{2}c_1 k^3}{w - 3\sqrt{2}k} \times \frac{\exp(\frac{6k^2 + w(w - 3\sqrt{2}k)}{k^2(w - 3\sqrt{2}k)}\xi)}{\frac{3c_1 k^4}{6k^2 + w(w - 3\sqrt{2}k)}\exp(\frac{6k^2 + w(w - 3\sqrt{2}k)}{k^2(w - 3\sqrt{2}k)}\xi) + c_2}, \qquad (22)$$

where $\xi = k(x \pm \frac{3}{\sqrt{2}}t)$ with $w = \pm\frac{3}{\sqrt{2}}k$. Here $c_1$ and $c_2$ are arbitrary constants.

If $c_2 = \frac{3k^4 c_1}{6k^2 + w(w - 3\sqrt{2}k)}$ then,

$$u = \pm 1 \pm \frac{\sqrt{2}\{6k^2 + w(w - 3\sqrt{2}k)\}}{k(w - 3\sqrt{2}k)}\left\{1 + \tanh\left(\frac{6k^2 + w(w - 3\sqrt{2}k)}{k^2(w - 3\sqrt{2}k)}\xi\right)\right\}, \qquad (23)$$

where $\xi = k(x \pm \frac{3}{\sqrt{2}}t)$ with $w = \pm\frac{3}{\sqrt{2}}k$.

If $c_2 = -\frac{3k^4 c_1}{6k^2 + w(w - 3\sqrt{2}k)}$ then,

$$u = \pm 1 \pm \frac{\sqrt{2}\{6k^2 + w(w - 3\sqrt{2}k)\}}{k(w - 3\sqrt{2}k)} \left\{1 + \tanh\left(\frac{6k^2 + w(w - 3\sqrt{2}k))}{k^2(w - 3\sqrt{2}k)}\xi\right)\right\}, \qquad (24)$$

where $\xi = k(x \pm \frac{3}{\sqrt{2}}t)$ with $w = \pm\frac{3}{\sqrt{2}}k$.

Since $c_1$ and $c_2$ are arbitrary constants, for other choices of $c_1$ and $c_2$ it might yield much new and more general exact solutions of the nonlinear Cahn Allen equation without any aid of symbolic computation software.

The solutions $u(x,t)$ obtained in Eqs. (23) is similar to the Fig-1 and (24) is similar to the Fig-2 and omitted for convenience.

Again with commercial software we can also find some solutions of the Cahn Allen equation (solving from (11) and (12)).

For $A_0 = 0$, $A_1 = \pm\sqrt{2}k$ we get $S(\xi) = a + b\exp(\pm\xi/\sqrt{2}k)$,

And thus $u(x,t) = \pm\dfrac{b}{a\left\{\cosh\dfrac{\xi}{\sqrt{2}k} \mp \sinh\dfrac{\xi}{\sqrt{2}k}\right\} + b}$ with $\xi = k\left(x \pm 3t/\sqrt{2}\right)$. (25)

If we consider $a/b = \exp(2c)$ then,

Eq.(25) reduces to well known solution $u(x,t) = \pm\dfrac{1}{2}\left\{1 + \tanh\left(\pm\dfrac{1}{\sqrt{2}}x + \dfrac{3}{2}t + c\right)\right\}$ (26)

For $A_0 = 1$, $A_1 = \pm\sqrt{2}k$ we get $S(\xi) = a + b\exp(\pm\xi/\sqrt{2}k)$,

And thus $u(x,t) = 1 - \dfrac{b}{a\left\{\cosh\dfrac{\xi}{\sqrt{2}k} \pm \sinh\dfrac{\xi}{\sqrt{2}k}\right\} + b}$ with $\xi = k\left(x \pm 3t/\sqrt{2}\right)$. (27)

If we consider $a/b = \exp(2c)$ then,

Eq.(26) reduces to well known solution $u(x,t) = \dfrac{1}{2}\left\{1 + \tanh\left(\pm\dfrac{1}{\sqrt{2}}x + \dfrac{3}{2}t + c\right)\right\}$ (28)

For $A_0 = -1$, $A_1 = \pm\sqrt{2}k$ we get $S(\xi) = a + b\exp(\mp\xi/\sqrt{2}k)$,

And thus $u(x,t) = -1 - \dfrac{b}{a\left\{\cosh\dfrac{\xi}{\sqrt{2}k} \mp \sinh\dfrac{\xi}{\sqrt{2}k}\right\} + b}$ with $\xi = k\left(x \mp 3t/\sqrt{2}\right)$. (29)

If we consider $a/b = \exp(2c)$ then,

Eq.(26) reduces to well known solution $u(x,t) = -\dfrac{1}{2}\left\{1 + \tanh\left(\pm\dfrac{1}{\sqrt{2}}x + \dfrac{3}{2}t + c\right)\right\}$ (30)

Since $a$ and $b$ are arbitrary constants, for other choices of $a$ and $b$ it might yield much new and more general exact solutions of the nonlinear Cahn Allen equation. When we choose $a/b = \exp(2c)$ get special type solution like Eq.(28) and Eq.(30), but for other choose $a$ and $b$ in different way we can get different type of solutions. Thus Eq.(28) and Eq.(30) are special type of our solutions.

Graph of the solutions Eq.(25), Eq.(27) and Eq.(29) represent kink type **Fig-1** for same positive/negative values of the arbitrary constants $c_1$ and $c_2$ like but single soliton like **Fig-2** for opposite values of them.

### 4. Comparison

In this section, we compare our solution with some well-known methods namely exp-function method and first integral method as follows:

**a) Comparison with Exp-method Ref. [30]:** Ugurlu (Ref. [30]) obtained some solutions of the Cahn Allen equation via exp-function method in which solutions $u_8, u_9$ are identical with our solutions Eq. (25) when $b=1, a=b_0$ and the other solutions are different with their solutions( For more see the Ref. [28]).

**b) Comparison with First Integral method Ref. [31]:** Tascan and Bekir (See Ref. [31]) obtained some solutions of the Cahn Allen equation via first integral method in which solutions Eq. (3.20) are identical with our solutions Eq. (25) (when in our study $a=b=1, k=-1/\sqrt{2}$ and in their study $c_0=0$) $u_8, u_9$ are identical with our solutions Eq. (25) when $b=1, a=b_0$ and the other solutions are different with their solutions. On the contrary by using the MSE method in this article we obtained four solutions with less calculation.

## 5. Conclusions

In this article, we have successfully implemented the MSE method to find the exact traveling wave solutions of the Cahn-Allen equation. Comparing the MSE method to other methods, we claim that the MSE method is straightforward, efficient, and can be used in many other nonlinear evolution equations. In the existing methods, such as, the (G'/G)-expansion method, the Exp-function method, the tanh-function method it is required to make use of the symbolic computation software, such as Mathematica or Maple to facilitate the complex algebraic computations. To solve non-linear evolution equations via MSE method no auxiliary equations is needed. On the other hand, via the MSE method, the exact and solitary wave solutions to these equations have been achieved without using any symbolic computation software because the method is very simple and has easy computations.